\documentclass{article}
\newtheorem{theorem}{Theorem}
\newtheorem{lemma}{Lemma}

\usepackage{graphicx}
\DeclareSymbolFont{AMSa}{U}{msa}{m}{n}
\DeclareMathDelimiter\ulcorner{\mathopen} {AMSa}{"70}{AMSa}{"70}
\DeclareMathDelimiter\urcorner{\mathclose}{AMSa}{"71}{AMSa}{"71}
\makeatletter
\def\uufill{$\m@th\mathopen\ulcorner\mkern-7mu%
  \cleaders\hbox{\rule[6pt]{1dd}{1dd}}\hfill
  \mkern-7mu\mathclose\urcorner$}
\def\overbrack#1{\vbox{\m@th\ialign{##\crcr
      \uufill\crcr\noalign{\kern-\p@\nointerlineskip}%
      $\hfil\displaystyle{#1}\hfil$\crcr}}}
\makeatother
\title{On a Proof of \\ the Goldbach Conjecture and the Twin Prime Conjecture}

\author{Sze Kui Ng
\\ {\small Department of Mathematics,
Hong Kong Baptist University, Hong Kong}\\{\small szekuing@yahoo.com.hk}}

\begin{document}
\date{}
\maketitle
\begin{abstract}
In this paper we use the connected sum operation on knots to show that there is a one-to-one relation between knots and numbers. In this relation prime knots are bijectively assigned with prime numbers such that the prime number $2$ corresponds to the trefoil knot.
From this relation we have a classification table of knots where knots are one-to-one assigned with numbers. 
Further this assignment for the $n$th induction step of the number $2^n$ is determined by this assignment for the previous $n-1$ steps. From this induction of assigning knots with numbers we can solve some problems in number theory such as the Goldbach Conjecture and the Twin Prime Conjecture.

\end{abstract}

\section{Introduction}\label{sec00}

It is well known that in number theory we have the classical result that each number can be uniquely factorized as a product of prime numbers \cite{Har1}-\cite{Sha}. On the other hand it is also well known that in knot theory we have the classical result that each knot can be uniquely factorized as the connected sum of prime knots \cite{Ada}-\cite{Rol}. In this paper we show that there is a deeper relation between these two factorizations that
we use 
the connected sum operation on knots to find out  
a one-to-one relation between knots and numbers. 
In this relation prime knots are bijectively assigned with prime numbers such that the prime number $2$ corresponds to the trefoil knot.
From this relation we have a classification table of knots where knots are one-to-one assigned with numbers.

In forming this relation between  knots and numbers
we show by induction on the number $n$ of $2^n$ that this assignment for the $n$th induction step (including the distribution of prime numbers in the $n$th induction step) is determined by this assignment for the previous $n-1$ steps.
We then use this induction of assigning knots with numbers  
to solve some problems in number theory such as the Goldbach Conjecture and the Twin Prime Conjecture. 

This paper is organized as follows. In section 2 and 3 we give a classification table of knots where each knot is assigned with an integer and that prime knots are bijectively assigned with prime numbers such that the prime number $2$ corresponds to the trefoil knot.
Further this assignment for the $n$th induction step of the number $2^n$ is determined by this assignment for the previous $n-1$ steps. 
By using this induction of assigning knots with numbers
we then in section 4 prove some conjectures in number theory such as the Goldbach Conjecture and the Twin Prime Conjecture.

\section{A Classification Table of Knots I}\label{sec1}

We shall use only the connected sum operation on knots to find out a relation between knots and numbers. 
For simplicity we use the positive integer $|m|$ to form a classification table of knots where $m$ is assigned to a knot while $-m$ is assigned to 
its mirror image if the knot is not equivalent to its mirror image.
Our main references on knots are \cite{Ada}-\cite{Rol}.

Let $\star$ denote the connected sum of two knots such that the
resulting total number of alternating crossings is equal to the
sum of alternating crossings of each of the two knots minus $2$.
As an example we have the reef knot (or the square knot) ${\bf
3_1\star 3_1}$ which is a composite knot composed with the knot
$\bf 3_1$ and its mirror image. 
This square knot has $6$ crossings
and $4$ alternating crossings. Then let $\times$ denote the
connected sum for two knots such that the resulting total number
of alternating crossings is equal to the sum of alternating
crossings of each of the two knots. As an example we have the granny
knot ${\bf 3_1\times 3_1}$ which is a composite knot composed with
two identical knots $\bf 3_1$. 
(For simplicity we use one notation $\bf 3_1$ to denote both the trefoil knot and its mirror image though these two knots are nonequivalent). This knot has $6$ alternating crossings which is equal to the total number of crossings. We have that the two operations $\star$ and $\times$ satisfy the commutative law and the associative law
 \cite{Ada}-\cite{Rol}. Further for each knot there is a unique factorization of this knot into a $\star$ and $\times$ operations of prime knots which is similar to the unique factorization of a number into a product of prime numbers \cite{Ada}-\cite{Rol}. In this paper we show that there is a deeper connection between these two factorizations.

 The first aim of this paper is to find out a table of the relation between knots and numbers by using only the operations $\star$ and $\times$ on knots and by using the following data as the initial step for induction:

{\bf Initial data for induction:} The prime knot ${\bf 3_1}$ is assigned with the number $1$ and it also plays the role of $2$. This means that the number $2$ is not assigned to other knots and is left for the prime knot ${\bf 3_1}$.
$\diamond$

Before the induction let us give remarks on this initial data for induction.

{\bf Remark}. The reader may wonder that why in this initial data the prime knot ${\bf 3_1}$ is not assigned directly with the prime number $2$. The reason for this initial data is that, by using a quantum gauge model approach for representing knots we show that the prime knot ${\bf 3_1}$ should be assigned with the number $1$ and it also plays the role of $2$ \cite{Ng}. This phenomeon is due to the property that the $\star$ operation is partially similar to the addition and partially similar to the multiplication where $1$ is for addition and $2$ is for multiplication.

In \cite{Ng} by using a quantum gauge model approach we also show that the next prime knot ${\bf 4_1}$ is assigned with the prime number $3$  and the composite knots 
${\bf 3_1\star 3_1}$ and ${\bf 3_1\times 3_1}$ are assigned with the number $4$ and $9$ respectively. This is as a guide for us to establish the classification table of knots. However in this paper we shall use only the elementary connected sum operation approach instead of the whole quantum gauge model approach to establish this classification table of knots.

{\bf Remark}. We shall say that the prime knot ${\bf 3_1}$ is assigned with the number $1$ and is related to the prime number $2$. 
$\diamond$ 

We shall give an induction on the number $n$ of $2^n$ for establishing the table. For each induction step on $n$ because of the special role of the trefoil knot ${\bf 3_1}$ we let the composite knot ${\bf 3_1}^n$ obtained by repeatedly taking $\star$ operation $n-1$ times on the trefoil knot ${\bf 3_1}$ be assigned with the number $2^n$ in this induction.

Let us first give the following table relating knots and numbers up to $2^5$ as a guide for the induction  
for establishing the whole classification table of knots:

\begin{displaymath}
\begin{array}{|c|c|c|c|} \hline
\mbox{Type of Knot}& \mbox{Assigned number} \,\, |m|
 &\mbox{Type of Knot}& \mbox{Assigned number} \,\, |m|
\\ \hline
{\bf 3_1} & 1 & {\bf 6_3} &  17\\ \hline

 &  2
&  {\bf 3_1\times 4_1} &  18 \\ \hline

{\bf 4_1} &  3 & {\bf 7_1} &  19 \\ \hline

{\bf 3_1\star 3_1} &  4 &
{\bf 4_1\star 5_1} &  20
\\ \hline

{\bf 5_1} & 5 & {\bf 4_1\star(3_1\star 4_1) } &  21
\\ \hline

{\bf 3_1\star 4_1} & 6 & {\bf 4_1\star 5_2} & 22 \\ \hline

{\bf 5_2} &  7 & {\bf 7_2} & 23 \\ \hline

{\bf 3_1\star 3_1\star 3_1} &  8 &
{\bf 3_1\star (3_1\times 3_1)}& 24 \\ \hline

{\bf 3_1\times 3_1} & 9 &
{\bf 3_1\star (3_1\star  5_1)}& 25 \\ \hline

{\bf 3_1\star 5_1} &  10 & {\bf 3_1\star 6_1}
& 26 \\ \hline

{\bf 6_1} &  11 & {\bf 3_1\star (3_1\star  5_2)}& 27
\\ \hline

{\bf 3_1\star 5_2} &  12 & {\bf 3_1\star 6_2} & 28
\\ \hline

{\bf 6_2} &  13 & {\bf 7_3} & 29 \\ \hline

{\bf 4_1\star 4_1} &  14 & { \bf (3_1\star 3_1)\star (3_1\star
4_1)} & 30 \\ \hline

{\bf 4_1\star (3_1\star 3_1)} & 15 & {\bf 7_4} & 31
\\ \hline

{ \bf (3_1\star 3_1)\star (3_1\star 3_1)} & 16 & {\bf (3_1\star
3_1)\star (3_1\star 3_1\star 3_1)} &  32
\\ \hline
\end{array}
\end{displaymath}

From this table we see that the $\star$ operation is similar to the usual multiplication $\cdot$ on numbers. Without the $\times$ operation this $\star$ operation would be exactly the usual multipilcation on numbers if this $\star$ operation is regarded as an operation on numbers.
From this table we see that comparable composite knots (in a sense
from the table and we shall discuss this point later) are grouped
in each of the intervals between two prime numbers. It is
interesting that in each interval composite numbers are one-to-one
assigned to the comparable composite knots while prime numbers are one-to-one
assigned to prime knots. Here a main point is to introduce the $\times$ operation while keeping composite knots correspond to composite numbers and prime knots correspond to prime numbers. To this end we need to have rooms at the positions of composite numbers for the introduction of composite knots obtained by the $\times$ operation. We shall show that these rooms can be obtained by using the special property of the trefoil knot which is assigned with the number $1$ (for the addition property of the $\star$ and $\times$ operations) while this trefoil knot is similar to the number $2$ for the multiplication property of the $\star$ operation.

Let us then carry out the induction steps for obtaining the whole table. To this end let us investigate in more detail the above comparable properties of knots. We have the following definitions and theorems.

{\bf Definition}.
We write $K_1<K_2$ if
$K_1$ is before $K_2$ in the ordering of knots; i.e. the number assigned to $K_1$ is less than the number assigned to $K_2$.

{\bf Definition (Preordering)}. Let two knots be
written in the form $K_1\star K_2$ and $K_1\star K_3$ where we
have determined the ordering of $K_2$ and $K_3$. Then we say that
$K_1\star K_2$ and $K_1\star K_3$ are in a preordering in the sense that we put the ordering of these two knots to follow the ordering of $K_2$ and $K_3$. If this preordering is not changed by conditions from other preorderings on these two knots (which are from other factorization forms of these two knots) then this preordering becomes the ordering of these two knots. 
We shall see that this preordering gives the comparable property in the above table. $ \diamond$

{\bf Remark}.
a) This definition is consistent since if $K_1$ is the unknot then we have 
$K_1\star K_2$=$K_2$ and $K_1\star K_3$=$K_3$ and thus the ordering of $K_1\star K_2$ and $K_1\star K_3$ follows the ordering of $K_2$ and $K_3$.

b) We can also define similarly the preordering of two knots $K_1\times K_2$ and 
$K_1\times K_3$ with the $\times$ operation. $\diamond$

We have the following theorem.
\begin{theorem}
Consider two knots of the form $K_1\star K_2$ and $K_1\star K_3$ where $K_1$, $K_2$ and $K_3$ are prime knots such that $K_2<K_3$. 
Then we have $K_1\star K_2 < K_1\star K_3$.
\label{pre}
\end{theorem}
{\bf Proof}.
Since $K_1$, $K_2$ and $K_3$ are prime knots there are no other factorization forms of the two knots $K_1\star K_2$ and $K_1\star K_3$. Thus these two forms of the two knots are the only way to give preordering to the two knots and thus there are no other conditions to change the preordering given by this factorization form of the two knots. Thus we have that $K_2<K_3$ implies $K_1\star K_2 < K_1\star K_3$.
$\diamond$

\begin{theorem}
Suppose two knots are written in the form $K_1\star K_2$ and $K_1\star K_3$ for determining their ordering and that the other forms of these two knots are not for determining their ordering. Suppose that $K_2<K_3$. Then we have $K_1\star K_2 < K_1\star K_3$.
\end{theorem}
{\bf Proof}.
The proof of this theorem is similar to the proof of the above theorem. Since the other factorization forms are not for the determination of the ordering of the two knots in the factorization form $K_1\star K_2$ and $K_1\star K_3$ we have that the preordering of these two knots in this factorization form becomes the ordering of these two knots. Thus we have $K_1\star K_2 < K_1\star K_3$. $\diamond$

As a generalization of theorem \ref{pre} we have the following theorem.

\begin{theorem}
Let two knots be of the form $K_1\star K_2$ and $K_1\star K_3$ where $K_2$ and $K_3$ are prime knots. Suppose that $K_2<K_3$. Then we have $K_1\star K_2 < K_1\star K_3$.
\label{pre2}
\end{theorem}
{\bf Proof}.
We have the preordering that $K_1\star K_2$ is before $K_1\star K_3$. Then since $K_2$ and $K_3$ are prime knots we have that the other preordering of $K_1\star K_2$ and $K_1\star K_3$ can only from the factorization of $K_1$. Without loss of generality let us suppose that $K_1$ is of the form $K_1=K_4\star K_5$ where 
$K_4< K_5$ and $K_4$ and $K_5$ are prime knots. Then we have the factorization 
$K_1\star K_2=K_4\star (K_5\star K_2)$ and $K_1\star K_3=K_5\star (K_4\star K_3)$. 
This factorization is the only factorization that might change the preordering that $K_1\star K_2$ is before $K_1\star K_3$. Then if $K_2\neq K_4$ or $K_3\neq K_5$
with this  factorization the two knots $K_1\star K_2$ and $K_1\star K_3$ are noncomparable in the sense that this factorization gives no preordering property and that the ordering of these two knots is determined by other conditions. Thus this factorization of $K_1\star K_2$ and $K_1\star K_3$ is not for the determination of the ordering of $K_1\star K_2$ and $K_1\star K_3$. Thus the preordering that $K_1\star K_2$ is before $K_1\star K_3$ is the ordering of $K_1\star K_2$ and $K_1\star K_3$. On the other hand if $K_2=K_4$ and $K_3=K_5$ then this factorization gives the same preordering that $K_1\star K_2$ is before $K_1\star K_3$. Thus for this case the preordering that $K_1\star K_2$ is before $K_1\star K_3$ is also the ordering of $K_1\star K_2$ and $K_1\star K_3$. Thus we have $K_1\star K_2 <K_1\star K_3$. $\diamond$

In addition to the above theorems we have the following theorems.
\begin{theorem}
Consider two knots of the form $K_1\times K_2$ and $K_1\times K_3$ where $K_1$, $K_2$ and $K_3$ are prime knots such that $K_2<K_3$. 
Then we have $K_1\times K_2 < K_1\times K_3$.
\end{theorem}
{\bf Proof}. By using a preordering property for knots with $\times$ operation as similar to that for knots with $\star$ operation we have that the proof of this theorem is similar to the proof of the above theorems. $\diamond$

\begin{theorem}
Let two knots be of the form $K_1\times K_2$ and $K_1\times K_3$ where $K_2$ and $K_3$ are prime knots. Suppose that $K_2<K_3$. Then we have $K_1\times K_2 < K_1\times K_3$.
\end{theorem}
{\bf Proof}. 
The proof of this theorem is also similar to the proof of the theorem \ref{pre2}. $\diamond$

These two theorems will be used for introducing and ordering knots involved with a $\times$ operation which will have the effect of pushing out composite knots with the property of jumping over (to be defined) such that knots are assigned with a prime number if and only if the knot is a prime knot.

Let us investigate more on the property of preordering. We consider the following

{\bf Definition (Preordering sequences)}
 At the $n$th induction step let the prime knot ${\bf 3_1}$ take a $\star$ operation with the previous $(n-1)$th step. We call this obtained sequence of composite knots as a
preordering sequence. Thus from the ordering of the $(n-1)$th step we have a sequence of composite knots which will be for the construction of the $n$th step. 

Then we let the prime knot ${\bf 4_1}$ (or the knot assigned with a prime number which is $3$ in the $2$nd step as can be seen from the above table) take a $\star$ operation with the previous $(n-2)$th step. From this we get a sequence of composite knots for constructing the $n$th step. 
Then we let the prime knots ${\bf 5_1}$ and ${\bf 5_2}$ (which are prime knots in the same step assigned with a prime number which is $5$ or $7$ in the $3$rd step as can be seen from the above table) take a $\star$ operation with the previous $(n-3)$th step respectively. From this we get two sequences for constructing the $n$th step.

Continuing in this way until the sequences are obtained by a prime knot in the $(n-1)$th step taking a $\star$ operation with the step $n=1$ where the prime knot is assigned with a prime number in the $(n-1)$th step by induction (By induction each prime number greater than $2$ will be assigned to a prime knot). 

We call these obtained sequences of composite knots as 
the preordering sequences of composite knots for constructing the $n$th step. Also we call
the sequences truncated from these preordering sequences as preordering subsequences of composite knots for constructing the $n$th step. $\diamond$

We first have the following lemma on preordering sequence.
\begin{lemma}
Let $K$ be a knot in a preordering sequence of the $n$th step. 
Then there exists a room for this $K$ in the $n$th step in the sense that this $K$ corresponds to a number in the $n$th step or in the $(n-1)$th step.
\end{lemma}
{\bf Proof}. Let $K$ be of the form $K= {\bf3_1}\star K_1$ where $K_1$ is a knot in the previous $(n-1)$th step. By induction we have that $K_1$ is assigned with a number $a$ which is the position of $K_1$ in the previous $(n-1)$th step. Then since ${\bf3_1}$ corresponds to the number $2$ we have that $K$ corresponds to the number $2\cdot a$ in the $n$th step (We remark that 
$K$ may not be assigned with the number $2\cdot a$). Thus there exists a room for this $K$ in the $n$th step.
Then let $K$ be of the form $K= {\bf4_1}\star K_2$ where $K_2$ is a knot in the previous $(n-2)$th step. By induction we have that $K_2$ is assigned with a number $b$ which is the position of $K_2$ in the previous $(n-2)$th step. Since ${\bf4_1}$ is by induction assigned with the prime number $3$ we have $3\cdot b >3\cdot2^{n-3}>2\cdot2^{n-3}=2^{n-2}$. Also we have 
$3\cdot b <3\cdot2^{n-2}<2^2\cdot2^{n-2}=2^{n}$. Thus there exists a room for this $K$ in the $(n-1)$th step or the $n$th step.
Continuing in this way we have that this lemma holds. $\diamond$

{\bf Remark}. By using  this lemma we shall construct each $n$th step of the classification table by first filling the $n$th step with the preordering subsequences of the $n$th step. $\diamond$

{\bf Remark}.
When the number corresponding to the knot $K$ in the above proof is not in the $n$th step we have that the knot $K$ in the preordering sequences of the $n$th step has the function of pushing a knot $K^{\prime}$ out of the $n$th step where this knot $K^{\prime}$ is related to a number in the $n$th step in order for the knot $K$ to be filled into the $n$th step. 

As an example in the above table the knot $K= {\bf4_1}\star{\bf5_1} $ (related to the number $3\cdot5$) in a preordering sequence of the $5$th step pushes the knot $K^{\prime}= {\bf5_1}\star {\bf5_1}$ related to the number $5\cdot5$ in the $5$th step out of the $5$th step. This relation of pushing out is by the chain $3\cdot5 \to 2\cdot2\cdot5 \to 5\cdot5$.

As another example in the above table the knot $K= {\bf3_1}\star({\bf3_1}\times{\bf3_1})$ (correspoded to the number $2\cdot9$) in a preordering sequence of the $5$th step pushes the knot $K^{\prime}= {\bf3_1}\star({\bf4_1}\star {\bf5_1})$ related to the number $2\cdot3\cdot5$ in the $5$th step out of the $5$th step. This relation of pushing out is by the chain $2\cdot9 \to 2\cdot2\cdot2\cdot3 \to 2\cdot3\cdot5$.
$\diamond$

\begin{lemma} 
For $n\geq 2$ the preordering subsequences for the $n$th step can cover the whole $n$th step.  
\end{lemma}
{\bf Proof}. 
For $n=2$ we have one preordering sequence with number of knots $=2^0$ which is obtained by the prime knot ${\bf 3_1}$ taking $\star$ operation with the step $n=2-1=1$.  
In addition we have the knot ${\bf 3_1}\star{\bf 3_1}$ which is assigned at the position of $2^n,n=2$ by the induction procedure.  Then since  the total rooms of this step $n=2$ is $2^1$ we have that
these two knots cover this step $n=2$.

For $n=3$ we have one preordering sequence with number of knots $=2^1$ which is obtained by the prime knot ${\bf 3_1}$ taking $\star$ operation with the step $3-1=2$. This sequence covers half of this step $n=3$ which is with $2^{3-1}=2^2$ rooms. 
Then we have one more preordering sequence  which is obtained by the knot ${\bf 4_1}$ taking $\star$ operation with step $n=1$ giving the number $2^0=1$ of knots. This covers half of the remaining rooms of the step $n=3$ which is with $2^{2-1}=2^1$ rooms. Then
in addition we have the knot ${\bf 3_1}\star{\bf 3_1}$ which is assigned at the position of $2^n,n=2$ by the induction procedure. These four knots (repeated knots are counted) thus cover the step $n=3$.

For the $n$th step we have one preordering sequence with the number of knots $=2^{n-2}$ which is obtained by the prime knot ${\bf 3_1}$ taking $\star$ operation with the $(n-1)$th step. This sequence covers half of this $n$th step which is with $2^{n-1}$ rooms. Then we have a preordering sequence  which is obtained by the knot ${\bf 4_1}$ taking $\star$ operation with the $(n-2)$th step giving the number $2^{n-3}$ of knots.  This covers half of the remaining rooms of the $n$th step which is with the remaining $2^{n-2}$ rooms. Then we have one preordering sequence obtained by picking a prime knot (e.g.${\bf 5_1}$) which by induction is assigned with a prime number (e.g. the number 5) taking $\star$ operation with the $(n-3)$th step.

Continuing in this way until the knot ${\bf 3_1}^n$ 
 is by induction assigned at the position of $2^n$. The total number of these knots is $2^{n-1}$ and thus these knots cover this $n$th step. This proves the lemma. $\diamond$

{\bf Remark}. 
Since there will have more than one prime numbers in the $k$th steps ($k>2$) 
in the covering  of the $n$th step there will have knots from 
the preordering sequences in repeat and in overlapping. These knots in repeat and in overlapping  may be deleted when  the ordering of the subsequences of the preordering sequences has been determinated for the covering of the $n$th step.

Also in the preordering sequences some knots which are in repeat and are not used for the covering of the $n$th step will be omitted when the ordering of the subsequences of the preordering sequences has been determinated for the covering of the $n$th step. $\diamond$

Let us then introduce another definition for constructing the classification table of knots.

{\bf Definition (Jumping over of the first kind)}. At an induction $n$th step consider a knot
$K^{\prime}$ and the knot $K={\bf 3_1}^n$ which is a $\star$ product of $n$ knots
$3_1$. $K^{\prime}$ is said to jump over $K$, denoted by $K \prec
K^{\prime}$,
 if exist $K_2$ and $K_3$  such that
$K^{\prime}=K_2\star K_3$ and for any $K_0$, $K_1$
such that  $K = K_0\star K_1$ where $K_0$, $K_1$, $K_2$ and $K_3$
are not
equal to ${\bf 3_1}$ we have
\begin{equation}
2^{n_0}<p_1\cdot\cdot\cdot p_{n_2}, \quad
2^{n_1}>q_1\cdot\cdot\cdot q_{n_3}
\label{class1}
\end{equation}
or vice versa
\begin{equation}
2^{n_0}>p_1\cdot\cdot\cdot p_{n_2}, \quad
2^{n_1}<q_1\cdot\cdot\cdot q_{n_3}
\label{class2}
\end{equation}
where $2^{n_0}$, $2^{n_1}$ are assigned to $K_0$ and $K_1$ respectively ($n_0 +n_1=n$) and
\begin{equation}
K_2= K_{p_1}\star \cdot\cdot\cdot \star K_{p_{n_2}} \quad K_3=
K_{q_1}\star \cdot\cdot\cdot \star K_{q_{n_3}}
\label{class3}
\end{equation}
where $K_{p_i}$, $K_{q_j}$ are prime knots which have been
assigned with prime integers $p_i$, $q_j$ respectively; and the following inequality holds:
\begin{equation}
2^n=2^{n_0+n_1}>p_1\cdot\cdot\cdot p_{n_2}\cdot q_1\cdot\cdot\cdot
q_{n_3}
\label{class33}
\end{equation}
Let us call this definition as the property of jumping over of the first kind. $\diamond$

We remark that the definition of jumping over of the first kind is a generalization of the above
ordering of ${\bf 4_1}\star{\bf 5_1}$ and ${\bf 3_1\star 3_1\star
3_1\star 3_1}$ in the above table in the step $n=4$ of $2^4$. 
Let us consider some examples of this
definition. Consider the knots $K^{\prime}=K_2\star K_3={\bf
4_1\star 5_1}$ and $K={\bf 3_1\star 3_1\star 3_1\star 3_1}$. For
any $K_0$, $K_1$ which are not equal to ${\bf 3_1}$ such that
$K=K_0\star K_1$ we have $2^{n_0}< 5$ and $2^{n_1}> 3$ (or vice
versa) where $3$, $5$ are the numbers of ${\bf 4_1}$ and
${\bf 5_1}$ respectively. Thus we have that ${\bf (3_1\star
3_1)\star(3_1\star 3_1)} \prec {\bf 4_1\star 5_1}$.

As another example we have that
$ {\bf 3_1\star(3_1\star 3_1)\star(3_1\star 3_1)} \prec
{\bf 5_1\star 5_1}$, ${\bf 4_1\star 4_1\star 4_1 }$,  and
${\bf 3_1\star (4_1\star 5_1)}$.

{\bf A Remark on Notation}. At the $n$th step let a composite knot of the form $K_1\star K_2\star\cdot\cdot\cdot\star K_q$ where each $K_i$ is a prime knot such that $K_i$ is assigned with a prime number $p_i$ in the previous $n-1$ steps. Then in general $K_1\star K_2\star\cdot\cdot\cdot\star K_q$ is not assigned with the number $p_1\cdot\cdot\cdot p_q$. However with a little confusion and for notation convenience we shall sometimes use the notation $p_1\cdot\cdot\cdot p_n$ to denote the knot $K_1\star K_2\star\cdot\cdot\cdot\star K_q$ and we say that this knot is related to the number $p_1\cdot\cdot\cdot p_n$ (as similar to the knot ${\bf3_1}$ which is related to the number $2$ but is assigned with the number $1$) and we keep in mind that the knot $K_1\star K_2\star\cdot\cdot\cdot\star K_q$ may not be assigned with the number $p_1\cdot\cdot\cdot p_n$. With this notation then we may say that the composite number $3\cdot5$ jumps over the number $2^4$ which means that the composite knot ${\bf4_1}\star{\bf5_1}$ jumps over the knot ${\bf3_1}\star{\bf3_1}\star{\bf3_1}\star{\bf3_1}$. $\diamond$

{\bf Definition (Jumping over of the general kind)}. At the $n$th step let a composite knot $K^{\prime}$ be related with a number $p_1\cdot p_2\cdot\cdot\cdot p_m$ where the number $p_1\cdot p_2\cdot\cdot\cdot p_m$ is in the $n$th step. Then we say that the knot $K^{\prime}$ (or the number $p_1\cdot p_2\cdot\cdot\cdot p_m$) is of jumping over of the general kind (with respect to the knot $K$ in the definition of the jumping over of the first kind and we also write $K\prec K^{\prime}$) if 
$K^{\prime}$ satisfies one of the following conditions:

1) $K^{\prime}$ (or the number related to $K^{\prime}$) is of jumping over of the first kind; or

2) There exists a $p_i$ (for simplicity let it be $p_1$) and a prime number $q$ such that $p_1$ and $q$ are in the same step $k$ for some $k$ and $q$ is the largest prime number in this step such that the numbers $p_1\cdot p_2\cdot\cdot\cdot p_m$ and $q\cdot p_2\cdot\cdot\cdot p_m$ are also in the same step and that the knot $K_q^{\prime}$ related with the number $q\cdot p_2\cdot\cdot\cdot p_m$ is of jumping over of the first kind. $\diamond$

{\bf Remark}. The condition 2) is a natural generalization of 1) that if $K^{\prime}$ and the knot $K_q^{\prime}$ are as in 2) then they are both in the preordering sequences of an induction $n$th step or both not. 
Then since $K_q^{\prime}$  is of jumping over into an $(n+1)$th induction step and thus is not in the preordering sequences of the induction $n$th step we have that $K^{\prime}$ is also of jumping over into this $(n+1)$th induction step (even if $K^{\prime}$ is not of jumping over of the first kind). This means that $K^{\prime}$ is of jumping over of the general kind. $\diamond$

{\bf Example of jumping over of the general kind}. At an induction step let $K^{\prime}$ be related with $11\cdot5\cdot5$ (where we let $p_1=11$) and let $K_q^{\prime}$ be related with $13\cdot5\cdot5$ (where we let $q=13$). Then $K_q^{\prime}$ is of jumping over of the first kind. Thus we have that $K^{\prime}$ is of jumping over (of the general kind). $\diamond$

We shall show that if $K={\bf 3_1}^n \prec K^{\prime}$ then we can set
$K={\bf 3_1}^n<K^{\prime}$. Thus we have, in the above first example, ${\bf
(3_1\star 3_1)\star(3_1\star 3_1)}< {\bf 4_1\star 5_1}$ while $2^4
> 3\cdot 5$.  From this
property we  shall have rooms for the introduction of the $\times$
operation such that composite numbers are assigned to composite knots and prime numbers are assigned to prime knots. We have the following theorem.

\begin{theorem}
If $K ={\bf 3_1}^n\prec K^{\prime}$ then it is consistent with the preordering property that $K={\bf 3_1}^n <K^{\prime}$ for establishing the table. 
\end{theorem}

For proving this theorem let us first prove the following lemma.
\begin{lemma}
The preordering sequences for the construction of the $n$th step do not have knots of jumping over of the general kind.
\end{lemma}
{\bf Proof of the lemma}.
It is clear that the preordering sequence obtained by the ${\bf 3_1}$ taking a $\star$ operation with the previous $(n-1)$th step has no knots with the jumping over of the first kind property since ${\bf 3_1}$ is corresponded with the number $2$ and the  previous $(n-1)$th step has no knots with the jumping over of the first kind property for this $(n-1)$th step. Then preordering sequence obtained by the ${\bf 4_1}$ taking a $\star$ operation with the previous $(n-2)$th step has no knots with the jump over of the first kind property since ${\bf 4_1}$ is assigned with the number $3$ and $3<2^2$ and the  previous $(n-2)$th step has no knots with the jumping over  of the first kind property for this $(n-2)$th step. Continuing in this way we have that all the knots in these preordering sequences do not satisfy the property of jumping over of the first kind. Then let us show that these preordering sequences have no knots with the property of jumping over of the general kind. Suppose this is not true. Then there exists a knot with the property of jumping over of the general kind and let this knot be related with a number of the form
$p_1\cdot p_2\cdot\cdot\cdot p_m$ as in the definition of jumping over of the general kind such that there exists a prime number $q$ and that $p_1$ and $q$ are in the same step $k$ for some $k$ and $q$ is the largest prime number in this step such that the numbers $p_1\cdot p_2\cdot\cdot\cdot p_m$ and $q\cdot p_2\cdot\cdot\cdot p_m$ are also in the same step and the knot $K_q$ represented by $q\cdot p_2\cdot\cdot\cdot p_m$ is of jumping over of the first kind. Then since $p_1$ and $q$ are in the same step $k$ such that the numbers $p_1\cdot p_2\cdot\cdot\cdot p_m$ and $q\cdot p_2\cdot\cdot\cdot p_m$ are also in the same step we have that the two knots related with $p_1\cdot p_2\cdot\cdot\cdot p_m$ and $q\cdot p_2\cdot\cdot\cdot p_m$ are elements of two preordering sequences for the construction of the same $n$th step. Now since we have shown that the preordering sequences for the construction of the $n$th step do not have knots of jumping over of the first kind we have that this is a contradiction. This proves the lemma.
$\diamond$

{\bf Proof of the theorem}.
By the above lemma if $K={\bf 3_1}^n\prec K^{\prime}$ then $ K^{\prime}$ is not in the preordering sequences for the $n$th step and thus is pushed out from the $n$th step by the preordering sequences for the $n$th step and thus we have $K={\bf 3_1}^n< K^{\prime}$, as was to be proved.
$\diamond$

{\bf Remark}. We remark that there may exist knots (or numbers related to the knots) which are not in the preordering sequences and are not of jumping over. An example of such special knot is the knot ${\bf4_1}\star{\bf5_1}\star{\bf5_1}$ related with $3\cdot 5\cdot 5$ (but is not assigned with this number). 
$\diamond$

{\bf Definition}.
When there exists a knot which is not in the preordering sequences of the $n$th step and is not of jumping over we put this knot back into the $n$th step to join the  preordering sequences for the filling and covering of the $n$th step. Let us call the preordering sequences together with the knots which are not in the preordering sequences of the $n$th step and are not of jumping over as the generalized preordering sequences (for the filling and covering of the $n$th step).
$\diamond$

{\bf Remark}. By using the generalized preordering sequences for the covering of the $n$th step we have that the knots (or the number related to the knots) in the $n$th step pushed out of the $n$th step by the generalized preordering sequences  are just the knots of jumping over (of the general kind). $\diamond$

Then we also have the following theorem.
\begin{theorem}
At each $n$th step ($n>3$) in the covering of the $n$th step ($n>3$) with the generalized preordering sequences there are rooms for introducing new knots with the $\times$ operations.
\label{times}
\end{theorem}
{\bf Proof}.
We want to show that at each $n$th step ($n>3$) there are rooms for
introducing new knots with the $\times$ operations.
At $n=4$ we have shown that there is the room at the position $9$ for introducing the knot ${\bf 3_1}\times{\bf 3_1}$ with the $\times$ operation.
Let us suppose that this property holds at an
induction step $n-1$. Let us then consider the induction step $n$.
For each $n$ because of the relation between $1$ and $2$ for ${\bf
3_1}$  as a part of the induction step $n$ the number $2^n$ is assigned to the knot
${\bf 3_1}^n$ which is a $\star$
product of $n$ ${\bf 3_1}$. Then we want to show that for this
induction step $n$ by using the $\prec$ property we have rooms for
introducing the $\times$ operation. Let $K^{\prime}$ be a knot
such that ${\bf 3_1}^{n-1}\prec K^{\prime}$
and $K^{\prime}=K_2\star K_3$
is as in the definition of $\prec$ of jumping over of the first kind 
such that $p_1\cdot\cdot\cdot p_{n_2}\cdot q_1\cdot\cdot\cdot q_{n_3}<2^{n-1}$ 
(e.g. for $n-1=4$
we have $K^4={\bf 3_1}\star{\bf 3_1}\star{\bf 3_1}\star{\bf 3_1}$
and $K^{\prime}=K_2\star K_3={\bf 4_1}\star{\bf 5_1}$). Then let
us consider $K^{\prime\prime}=({\bf 3_1}\star K_2)\star K_3$.
Clearly we have ${\bf 3_1}^n\prec K^{\prime\prime}$. Thus for each
$K^{\prime}$  we have a $K^{\prime\prime}$ such that
${\bf 3_1}^n\prec K^{\prime\prime}$. Clearly all these $ K^{\prime\prime}$ are different.

Then from $K^{\prime}$ let us construct one more $K^{\prime\prime}$, as follows. Let $K^{\prime}$ be a knot of jumping over of the first kind. Let
$p_1\cdot\cdot\cdot p_{n_2}$ and $q_1\cdot\cdot\cdot q_{n_3}$ be
as in the definition of jumping over
of the first kind. Then as in the definition of jumping over
of the first kind (w.l.o.g)
we let 
\begin{equation}
2^{n_0}<p_1\cdot\cdot\cdot p_{n_2} \quad \mbox {and}\quad 
2^{n_1}>q_1\cdot\cdot\cdot q_{n_3}
\label{forjumpingover}
\end{equation} 

Then we have
\begin{equation}
2^{n_0+1}<(2\cdot p_1\cdot\cdot\cdot p_{n_2})-1\quad \mbox {and}\quad 
2^{n_1}>q_1\cdot\cdot\cdot q_{n_3}
\label{jumpover12}
\end{equation}
Also it is trivial that we have
$2^{n_0}<(2\cdot p_1\cdot\cdot\cdot p_{n_2})-1$ and
$2^{n_1+1}>q_1\cdot\cdot\cdot q_{n_3}$. This shows that ${\bf 3_1}^n\prec
K^{\prime\prime}:= K_{2a}\star K_{3}$ where $K_{2a}$ denotes the
knot with the number $(2\cdot p_1\cdot\cdot\cdot p_{n_2})-1$ as in
the definition of jumping over of the first kind (We remark that this $K^{\prime\prime}$
corresponds to the knot ${\bf 4_1}\star({\bf 4_1}\star{\bf 4_1})$
in the above induction step where $K_{2a}={\bf 4_1}\star{\bf 4_1}$
is with the number $2\cdot 5-1=3\cdot 3$). 
It is clear that all these more $K^{\prime\prime}$ are different from the above $K^{\prime\prime}$ constructed by the above method of taking a $\star$ operation with ${\bf 3_1}$.  
Thus there are more $ K^{\prime\prime}$ than $K^{\prime}$. Thus at this $n$th step there are rooms for introducing new knots with the $\times$ operations. 
This proves the theorem. $\diamond$

{\bf Remark}.
In the proof of the above theorem we have a way to construct the knots $ K^{\prime\prime}$ by replacing a number $a$ with the number $2a-1$. There is another way of constructing 
the knots $ K^{\prime\prime}$ by replacing a number $b$ with the number $2b+1$. For this way we need to check that the number related to $ K^{\prime\prime}$ is in the $(n-1)$th step for $ K^{\prime\prime}$ of jumping over into the $n$th step.
As an example let us consider the knot $ K^{\prime}={\bf 4_1}\star{\bf 4_1}\star{\bf 4_1}$ of jumping over into the $6$th step with the following data: 
\begin{equation}
2^{3}<3\cdot3 \quad \mbox {and}\quad 
2^{2}>3
\label{jumpingover9}
\end{equation}

From this data we have:
\begin{equation}
2^{3+1}<2\cdot3\cdot3-1 =17\quad \mbox {and}\quad 
2^{2}>3
\label{jumpingover10a}
\end{equation}
This data gives a knot $ K^{\prime\prime}$ with the related number $3\cdot17 $.
On the other hand from the data (\ref{jumpingover9}) we have:
\begin{equation}
2^{3}<3\cdot3 \quad \mbox {and}\quad 
2^{2+1}>2\cdot3+1
\label{jumpingover10}
\end{equation}
Since $(3\cdot3)(2\cdot3+1)=(2\cdot5-1)(2\cdot3+1)=2\cdot5\cdot2\cdot3+2\cdot2-1<2\cdot2\cdot2^4-1<2^6$ we have that the knot $ K^{\prime\prime}={\bf 4_1}\star{\bf 4_1}\star{\bf 5_2}$ related with the number $3\cdot3\cdot7$ is of jumping over into the $7$th step (We shall show that ${\bf 5_2}$
is assigned with the number $7$). $\diamond$

{\bf Remark}. The above theorem shows that at each $n$th step there are rooms for introducing new knots with the $\times$ operations and thus we may establish a one-to-one correspondence of knots and numbers such that prime knots are bijectively assigned with prime numbers. Further to this theorem we have the following main theorem:

\begin{theorem}
A classification table of knots can be formed (as partly described by the above table up to $2^n$ with $n=5$) by induction on the number $2^n$
such that knots are one-to-one assigned with an integer and prime knots are bijectively assigned with prime numbers such that the prime number $2$ corresponds to the trefoil knot. This assignment is onto the set of positive integers except $2$ where the trefoil knot is assigned with 1 and is related to $2$ and at each $n$th induction step of the number $2^n$ there are rooms for introducing new knots with the $\times$ operations only.

Further this assignment of knots to numbers for the $n$th induction step of the number $2^n$ effectively includes the determination of the distribution of prime numbers in the $n$th induction step and is by induction determined by this assignment for the previous $n-1$ induction steps such that the assignment for the previous $n-1$ induction steps is inherited in this assignment for the $n$th induction step as the preordering sequences in the determination of this assignment for the $n$th induction step.
\label{maintheorem}
\end{theorem}

{\bf Remark}. Let us also call this assignment of knots to numbers as the structure of numbers obtained by assigning numbers to knots. This structure of numbers is the original number system together with the one-to-one assignment of numbers to knots.
 
{\bf Proof}.
By the above lemmas and theorems we have that
the generalized preordering sequences 
have the function of pushing out those composite knots 
of jumping over from the $n$th step.
It follows that for step $n>3$ there must exist
%(and will also be unique by the preordering property)
 chains of transitions whose initial states are composite knots in repeat (to be replaced by the new composite knots with the $\times$ operations only); or the knots  
 jumping over into this $n$th step from the previous $(n-1)$th step or the knots in the preordering sequences with the $\times$ operations; such that the composite knots 
% (or the composite numbers related to these knots) 
 of jumping over are pushed out from the $n$th step by these chains. These chains are obtained by ordering the subsequences of the preordering sequences
 such that the preordering property holds in the %whole 
 $n$th step. 
 Further 
the intermediate states of the chains must be positions of composite numbers. This is because that if a chain is transited to an intermediate state which is a position of prime number then there are no composite knots related with  
this prime number and thus this chain can not be transited to the next state and  is stayed at the intermediate state forever and thus the chain can not push out the composite knot of jumping over. Then when a composite knot is transited to the position of an intermediate state (which is a position of composite number as has just been proved) this knot is definitely assigned with this composite number.
Then when a composite knot which is in repeat is transited to the position of an intermediate state this knot is also definitely assigned with this composite number. 
It follows that when the chains are completed we have that the ordering of the subsequences of preordering sequences is determined.

Then the remaining knots (which are not at the transition states of the chains) which are not in repeat are definitely assigned with the number of the position of these knots in the $n$th step. For these knots the numbers of positions assigned to them are just the number related to them respectively.

Then the remaining knots (which are not at the transition states of the chains) which are in repeat must be replaced by new prime knots because of the repeat and that no other composite knots related with numbers in this $n$th step in the generalized preordering sequences can be used to replace the remaining knots.  
This means that the numbers of the positions of these remaining knots in repeat are prime numbers in this $n$th step. This is because that if the number of the position assigned to the new prime knot is a composite number then the composite knot related  
with this composite number is either in a transition state or is not in transition.
If the composite knot is not in transition then the composite number related to 
 this composite knot is just the number assigning to this composite knot and since this number is  also assigned to the new prime knot that this is a contradiction. Then if this composite knot is in transition state then this means that the remaining knot is also in transition state and this is a contradiction since by definition the remaining knot is not at the transition states of the chains. 
 
 Thus prime numbers in the $n$th step are assigned and  are only assigned to prime knots which replace the remaining knots in repeat 
 in the $n$th step. 
 Thus from the preordering sequences we have determined the positions (i.e. the distribution) of prime numbers in the $n$th step. Now since the preordering sequences are constructed by  the previous steps  we have shown that the basic structure (in the sense of above proof) of this assignment of knots with numbers for the $n$th step (including the determination of the distribution of prime numbers in the $n$th step) is determined by this assignment of knots with numbers for the previous $n-1$ steps. In other words we have that the basic structure of the $n$th induction step  
 is determined by the structure of the previous $n-1$ steps.
 
 To complete the proof of this theorem 
 let us  show that at each $n$th induction step ($n>3$) there are rooms for introducing new composite knots with the $\times$ operations only and we can determine the ordering of these composite knots with the $\times$ operations only in each $n$th induction step. 
 
In the above proof we have shown that the basic structure of the $n$th induction step is determined by the structure of the previous steps such that the positions of the composite knots with the $\times$ operations only in the $n$th induction step are 
determined by the structures of the previous steps. 
These positions are fitted for the corrected composite knots with the $\times$ operations only constructed (by the $\times$ operations) by knots in the previous steps.  Thus for this $n$th induction step the introducing and ordering of composite knots with the $\times$ operations only is also determined by the structures of the previous $n-1$ steps.
 
Further since
the structures of the previous steps are inherited in the structure of the $n$th induction step as the preordering sequences in the determination of the structure of the $n$th induction step we have that all the properties of the structures of the previous steps are inherited in the structure of the $n$th induction step in the determination of the structure of the $n$th induction step. Thus the new composite knots with the $\times$ operations only in the $n$th induction step 
inherit the ordering properties (such as the preordering property) of composite knots with the $\times$ operations only  in the previous steps. 
(These ordering properties of the composite knots with the $\times$ operations only can be used to find out the corrected composite knots with the $\times$ operations only to be assigned at the corrected positions in the $n$th step).

With this fact let  
us then show that at each $n$th induction step ($n>3$) there are rooms for introducing new composite knots with the $\times$ operations only.
 As in the proof of the theorem \ref{times} we first construct more $ K^{\prime\prime}$ by the method following (\ref{forjumpingover}). Let us start at the step $n=4$. For this step we have the knot $K^{\prime}={\bf 4_1}\star{\bf 5_1}$ jumps over into the step $n=5$. For this $K^{\prime}$ we have the following data as in (\ref{forjumpingover}):
\begin{equation}
2^{2}<5 \quad \mbox {and}\quad 
2^{2}>3
\label{jumpingover3}
\end{equation}
From (\ref{jumpingover3}) we construct a $ K^{\prime\prime}$ for the step $n=5$ by the following data:
\begin{equation}
2^{2+1}<2\cdot5-1= 3\cdot3 \quad \mbox {and}\quad 
2^{2}>3
\label{jumpingover4}
\end{equation}
This data gives one more $ K^{\prime\prime}={\bf 4_1}\star{\bf 4_1}\star{\bf 4_1}$. Then from (\ref{jumpingover3}) we construct one more $ K^{\prime\prime}$ for the step $n=5$ by the following data:
\begin{equation}
2^{3}>5 \quad \mbox {and}\quad 
2^{1+1}<2\cdot3-1=5
\label{jumpingover5}
\end{equation}
This data gives one more $K^{\prime\prime}={\bf 5_1}\star{\bf 5_1}$. Thus in this step $n=5$ there are two rooms  for the two knots $ K^{\prime}={\bf 4_1}\star{\bf 5_1}$ and ${\bf 3_1}\star({\bf 3_1}\times{\bf 3_1})$ coming from the preordering sequences  and there exists exactly one room for introducing a  new composite knot with the $\times$ operations only (Recall that we also have a $K^{\prime\prime}={\bf 3_1}\star{\bf 4_1}\star{\bf 5_1}$). 
From the ordering of knots in the previous steps we determine that ${\bf 3_1}\times{\bf 4_1}$ is the composite knot with the $\times$ operations only for this step.
Thus at the $4$th and $5$th steps we can and only can introduce exactly one composite knot with the $\times$ operations only and they are the knots ${\bf 3_1}\times {\bf 3_1}$ and ${\bf 3_1}\times {\bf 4_1}$ respectively.
This shows that at the $4$th and the $5$th steps we can determine the number of prime knots with the minimal number of crossings $=3$ and $=4$ respectively (These two prime knots are denoted by ${\bf 3_1}$ and ${\bf 4_1}$ respectively and we do not distinguish knots with their mirror images for this determination of the ordering of knots with the $\times$ operations only. This also shows that there are rooms for introducing new composite knots with the $\times$ operations only in the $4$th and $5$th steps).

Then since this property is inherited
in the $6$th step  
we can thus determine that the $6$th step is a step for introducing new composite knots with the $\times$ operations only of the form ${\bf 3_1}\times {\bf 5_{(\cdot)}}$ where ${\bf 5_{(\cdot)}}$ denotes a prime knot with the minimal number of crossings $=5$ (and thus there are rooms for introducing new composite knots with the $\times$ operations only in this $6$th step). Also since the properties in the $ 4$th and $5$th steps are inherited 
in the $6$th step we can determine
the number of prime knots with the minimal number of crossings $=5$ by the knots of the form ${\bf 3_1}\times {\bf 5_{(\cdot)}}$ as
this is a property of knots with the $\times$ operations only in the $4$th and $5$th steps (In the classification table  in the next section we show that there are exactly two composite knots of the form ${\bf 3_1}\times {\bf 5_1}$ and ${\bf 3_1}\times {\bf 5_2}$ in the $6$th step whose ordering are determined by the preordering property of knots and the structure of the $6$th step. This thus shows that there are exactly two prime knots with
the minimal number of crossings $=5$ and they are denoted by ${\bf 5_1}$ and ${\bf 5_2}$ respectively).

Then since the properties of the $4$th, $5$th and $6$th steps are inherited %preserved 
in the  $7$th step we can determine that  the $7$th step is a step for introducing new composite knots with the $\times$ operations only of the form ${\bf 3_1}\times {\bf 6_{(\cdot)}}$ where ${\bf 6_{(\cdot)}}$ denotes a prime knot with the minimal number of crossings $=6$ (and thus there are rooms for introducing new composite knots with the $\times$ operations only in this $7$th step). Also since the properties in the $4$th, $5$th and $6$th steps are inherited %preserved 
in the $7$th step we can determine
the number of prime knots with the minimal number of crossings $=6$ by the knots of the form ${\bf 3_1}\times {\bf 6_{(\cdot)}}$ as 
this is a property of knots with the $\times$ operations only in the $4$th, $5$th and $6$th steps (In the classification table  in the next section we show that there are exactly three composite knots of the form ${\bf 3_1}\times {\bf 6_1}$, ${\bf 3_1}\times {\bf 6_2}$ and ${\bf 3_1}\times {\bf 6_3}$ in the $7$th step whose ordering are determined by the preordering property of knots and the structure of the $7$th step. This thus shows that there are exactly three prime knots with
the minimal number of crossings $=6$ and they are denoted by ${\bf 6_1}$, ${\bf 6_2}$ and ${\bf 6_3}$ respectively).

Continuing in this way we thus have that at each $n$th induction step $(n>3)$ we can determine the number of prime knots with the minimal number of crossings $=n-1$ and there are rooms for introducing new composite knots with the $\times$ operations only. This proves the theorem. $\diamond$

%On the other hand we have the ordering property that composite knots with the $\times$ operations only can be introduced and ordered in the previous $n-1$ steps and this ordering property (including the preordering property) is inherited in the $n$th induction step. Thus composite knots with the $\times$ operations only can be introduced and ordered in the $n$th induction step.
 
%This proves the theorem. $\diamond$

{\bf Example}. Let us consider the above table up to $2^5$ (with $n$ up to $5$) as an example.
For the induction step $n=2$ (or the induction step $2^2$ where we use $n=2$ to mean $2^n$ with $n=2$) we have one preordering sequence obtained by letting ${\bf 3_1}$ to take a $\star$ operation with the step $n=1$ (For the step $n=1$ the number $2^1$ is related to the trefoil knot ${\bf 3_1}$): ${\bf 3_1}\star {\bf 3_1}$. Then we fill the step $n=2$ with this preordering sequence and we have the following ordering of knots for this step $n=2$:
\begin{equation}
{\bf 3_1}\star {\bf 3_1},
{\bf 3_1}\star  {\bf 3_1}
\label{step2}
\end{equation}
where the first ${\bf 3_1}\star  {\bf 3_1}$ placed at the position $3$ is the preordering sequence while the second ${\bf 3_1}\star  {\bf 3_1}$ placed at the position $2^2$ is required by the induction procedure. For this step there is no numbers of jumping over.  Then we have that the first ${\bf 3_1}\star  {\bf 3_1}$ is a repeat of the second ${\bf 3_1}\star  {\bf 3_1}$.  
Thus this repeat one must be replaced by a new prime knot. Let us choose the prime knot ${\bf 4_1}$ to be this new prime knot since ${\bf 4_1}$ is the smallest of prime knots other than the trefoil knot. Then this new prime knot must be at the position of a prime number, as we have proved in the above theorem. Thus we have determined that $3$ is a prime number in this step $n=2$ by using the structure of numbers of step $n=1$ which is only with the prime number $2$.

Then for the induction step $n=3$ (or the induction step $2^3$) 
we have two preordering sequence obtained by letting ${\bf 4_1}$ to take a $\star$ operation with the step $n=1$ and by letting ${\bf 3_1}$ to take a $\star$ operation with the step $n=2$: 
\begin{equation}
{\bf 4_1}\star{\bf 3_1}; {\bf 3_1}\star{\bf 4_1},
{\bf 3_1}\star ({\bf 3_1}\star {\bf 3_1})
\label{step3a}
\end{equation}
where the first knot is the preordering sequence obtained by letting ${\bf 4_1}$ to take a $\star$ operation with the step $n=1$ and the second and third knots is the preordering sequence obtained by letting ${\bf 3_1}$ to take a $\star$ operation with the step $n=2$. 
For this step there is no numbers of jumping over and thus there are no chains of transition. Thus the ordering of the above three knots in this step follow the usual ordering of numbers.  
Thus the number assigned to the knot ${\bf 4_1}\star {\bf 3_1}={\bf 3_1}\star{\bf 4_1}$ must be assigned with a number less than that of ${\bf 3_1}\star {\bf 3_1}\star {\bf 3_1}$ by the ordering of ${\bf 3_1}\star{\bf 4_1}$ and ${\bf 3_1}\star {\bf 3_1}\star {\bf 3_1}$ in the second preordering sequence. 
By this ordering of the two preordering sequences  we have that the step $n=3$ is of the following form:
\begin{equation}
{\bf 4_1}\star{\bf 3_1};{\bf 3_1}\star{\bf 4_1},
{\bf 3_1}\star ({\bf 3_1}\star {\bf 3_1});
{\bf 3_1}\star {\bf 3_1}\star {\bf 3_1}
\label{cover}
\end{equation}
where the fourth knot ${\bf 3_1}\star {\bf 3_1}\star {\bf 3_1}$ is put at the position of $2^3$ and is assigned with the number $2^3$ as required by the induction procedure.
 Thus the third knot ${\bf 3_1}\star ({\bf 3_1}\star {\bf 3_1})$  is a repeated one and thus must be replaced by a prime knot and the position of this prime knot is determined to be a prime number.  Thus we have determined that the number $7$ is a prime number.   
 Then since there are no chains of transition we have that the composite knot ${\bf 3_1}\star{\bf 4_1}$ must be assigned with the number related to this knot and this number is $2\cdot 3=6$. Thus the composite knot ${\bf 3_1}\star{\bf 4_1}$ is at the position of $6$ and that 
the first knot ${\bf 4_1}\star{\bf 3_1}$ is a repeat of the second knot and thus must be replaced by a prime knot. Then since this prime knot is at the position of $5$ we have that $5$ is determined to be a prime number. Now the two prime knots at $5$ and $7$ must be the prime knots ${\bf 5_1}$ and ${\bf 5_2}$ respectively since these two knots are  the smallest prime knots other than ${\bf 3_1}$ and ${\bf 4_1}$ (We may just put in two prime knots and then later determine what these two knots will be. If we put in other prime knots then this will not change the distribution of prime numbers determined by the structure of numbers of the previous steps 
and it is only that the prime knots are assigned with incorrect prime numbers. Further as shown in the above proof by using knots of the form ${\bf 3_1}\times{\bf 5_{(\cdot)}}$  we can determine that there are exactly two prime knots with minimal number of crossings $=5$ and they are denoted by ${\bf 5_1}$ and ${\bf 5_2}$ respectively. From this we can then determine that these two prime knots are ${\bf 5_1}$ and ${\bf 5_2}$).
 Thus we have
the following ordering for $n=3$:
\begin{equation} 
{\bf 5_1}<{\bf 3_1}\star{\bf 4_1}<
{\bf 5_2}<{\bf 3_1}\star {\bf 3_1}\star {\bf 3_1}
\label{order1}
\end{equation} 
where ${\bf 5_1}$ is assigned with the prime number $5$ and ${\bf 5_2}$ is assigned with the prime number $7$. This gives the induction step $n=3$. For this step there is no knot with $\times$ operation since there is no knots of jumping over.

Let us then consider the step $n=4$ (or $2^4$). For this step we have the following three preordering sequences obtained from the steps $n=1,2,3$: 
\begin{equation}
\begin{array}{rl}
&{\bf 5_1}\star {\bf3_1};\\
&{\bf 4_1}\star {\bf 4_1},{\bf 4_1}\star {\bf 3_1}\star {\bf 3_1};\\
&{\bf 3_1}\star {\bf5_1},{\bf 3_1}\star {\bf 3_1}\star {\bf 4_1},{\bf 3_1}\star {\bf 5_2},{\bf 3_1}\star {\bf 3_1}\star {\bf 3_1}\star{\bf 3_1};\\
 \end{array}
\label{order2}
\end{equation}
where the third sequence is obtained by taking $\star$ operation of the knot ${\bf 3_1}$ with step $n=3$ while the second sequence is obtained by taking $\star$ operation of the knot ${\bf 4_1}$ with the step $n=2$ and the first sequence is obtained by taking $\star$ operation of the knot ${\bf 5_1}$ with step $n=1$. Then as required by the induction procedure the knot ${\bf 3_1}\star {\bf 3_1}\star {\bf 3_1}\star {\bf 3_1}$ is assigned at the position of $2^4$. The total number of knots in (\ref{order2}) plus this knot is exactly $2^3$ which is the total number of this step $n=4$. 

{\bf Remark}. We  have one more preordering sequence  obtained by taking $\star$ operation of the knot ${\bf 5_2}$ with step $n=1$. This preordering sequence gives the knot
${\bf 5_2}\star {\bf3_1}$. However since the knots in (\ref{order2}) and the knot ${\bf 3_1}\star {\bf 3_1}\star {\bf 3_1}\star {\bf 3_1}$ assigned at the position of $2^4$ are enough for covering this step $n=4$ and that the knot ${\bf 5_2}\star {\bf3_1}$ of this preordering sequence is a repeat of the knot ${\bf 5_2}\star {\bf3_1}$ in (\ref{order2}) that this preordering sequence obtained by taking $\star$ operation of the knot ${\bf 5_2}$ with step $n=1$ can be omitted. $\diamond$

Then to find the chains of transition for this step let us order the three preordering sequences with the following ordering where we rewrite the preordering sequences in column form and the knot ${\bf 3_1}\star {\bf 3_1}\star {\bf 3_1}\star {\bf 3_1}$ assigned at the position of $2^4$ is put to follow the three sequences:
\begin{equation}
\begin{array}{rl}
&{\bf 5_1}\star {\bf3_1};\\
&{\bf 3_1}\star {\bf5_1},\\
&{\bf 3_1}\star {\bf 3_1}\star {\bf 4_1},\\
&{\bf 3_1}\star {\bf 5_2},\\
&{\bf 3_1}\star {\bf 3_1}\star {\bf 3_1}\star{\bf 3_1};\\
&{\bf 4_1}\star {\bf 4_1},\\
&{\bf 4_1}\star {\bf 3_1}\star {\bf 3_1};\\
&{\bf 3_1}\star {\bf 3_1}\star {\bf 3_1}\star{\bf 3_1}
 \end{array}
\label{step4}
\end{equation}
We notice that this column exactly fills the step $n=4$.

For this step we have that the number $3\cdot5$ (or the knot ${\bf 4_1}\star {\bf 5_1}$ related with $3\cdot5$ ) is of jumping over.
From (\ref{step4}) we have the following chain of transition for pushing out ${\bf 4_1}\star {\bf 5_1}$ at $3\cdot5$ by a knot with the $\times$ operation replacing the repeated knot ${\bf 5_1}\star {\bf 3_1}$ at the position of $9=3\cdot3$:
\begin{equation}
{\bf 3_1}\times {\bf 3_1} (\mbox {at} 3\cdot3)\to {\bf 4_1}\star {\bf 4_1}(\mbox {at} 2\cdot7) \to{\bf 3_1}\star {\bf 5_2} (\mbox {at} 2\cdot2\cdot3)\to
{\bf 3_1}\star {\bf 3_1}\star {\bf 4_1} (\mbox {at} 3\cdot5)\to {\bf 4_1}\star {\bf 5_1} (\mbox {pushed out})
\label{Chain2}
\end{equation}
where we choose the knot ${\bf 3_1}\times{\bf 3_1}$ as the knot with the $\times$ operation since ${\bf 3_1}\times{\bf 3_1}$ is the smallest one of such knots.
For this chain the intermediate states are at positions of composite numbers $2\cdot7$, $2\cdot2\cdot3$ and $3\cdot5$. Thus the knots in this chain at the positions of these composite numbers are assigned with these composite numbers respectively.

Then once this chain of pushing out ${\bf 4_1}\star {\bf 5_1}$ at $3\cdot5$ is established we have that the other knots in repeat must by replaced by prime knots and that their positions must be prime numbers.
These positions are at $11$ and $13$ and thus $11$ and $13$ are determined to be prime numbers (The knot ${\bf 3_1}\star{\bf 3_1}\star{\bf 3_1}\star{\bf 3_1}$ at the end of this step must be assigned with $2^4=16$ by the induction procedure and thus the knot at $13$ is a repeat).
Then  the new prime knots ${\bf 6_1}$ and ${\bf 6_2}$ are suitable knots corresponding to the prime numbers $11$ and $13$ respectively since they are the smallest  prime knots other than ${\bf 3_1}$. ${\bf 4_1}$, ${\bf 5_1}$ and ${\bf 5_2}$ (As the above induction step we may just put in two prime knots and then later determine what these two prime knots will be. As shown in the above proof by using knots of the form 
${\bf 3_1}\times{\bf 6_{(\cdot)}}$ we can determine that there are exactly three prime knots with minimal number of crossings $=6$ and they are denoted by ${\bf 6_1}$, ${\bf 6_2}$ and ${\bf 6_3}$ respectively. From this we can then determine that these two prime knots are ${\bf 6_1}$ and ${\bf 6_2}$).

 This completes the step $n=4$.  Thus the structure of numbers 
of this step (including distribution of prime numbers in this step) is  determined by the structure of numbers of the previous induction steps.

Let us then consider the step $n=5$. For this step we have the following four preordering sequences from the previous steps $n=1,2,3,4$:
\begin{equation}
{\bf 6_1}\star{\bf 3_1}
\label{step5a}
\end{equation}
and
\begin{equation}
\begin{array}{rl}
 & {\bf 5_2}\star{\bf 4_1},\\
 & {\bf 5_2}\star ({\bf 3_1}\star {\bf 3_1})
\end{array}
\label{step5b}
\end{equation}
and
\begin{equation}
\begin{array}{rl}
 & {\bf 4_1}\star{\bf 5_1},\\
 & {\bf 4_1}\star ({\bf 3_1}\star{\bf 4_1}),\\
 & {\bf 4_1}\star {\bf 5_2},\\
 & {\bf 4_1}\star ({\bf 3_1}\star {\bf 3_1}\star {\bf 3_1})
 \end{array}
\label{CC2}
\end{equation} 
and
\begin{equation}
\begin{array}{rl}
 & {\bf 3_1}\star ({\bf 3_1}\times {\bf 3_1}),\\
 & {\bf 3_1}\star({\bf 3_1}\star {\bf 5_1}),\\
 & {\bf 3_1}\star {\bf 6_1},\\
 & {\bf 3_1}\star ({\bf 3_1}\star{\bf 5_2}),\\
 & {\bf 3_1}\star{\bf 6_2}, \\
 & {\bf 3_1}\star ({\bf 4_1}\star {\bf 4_1}),\\ 
 & {\bf 3_1}\star ({\bf 3_1}\star {\bf 3_1}\star {\bf 4_1}),\\
 & {\bf 3_1}\star ({\bf 3_1}\star {\bf 3_1}\star {\bf 3_1}\star {\bf 3_1})
\end{array}
\label{CC4}
\end{equation}

The total number of knots (including repeat) in the above sequences plus the knot ${\bf 3_1}\star {\bf 3_1}\star {\bf 3_1}\star {\bf 3_1}\star {\bf 3_1}$ to be assigned at the position of $2^5$ exactly cover this $n=5$ step. 

{\bf Remark}. As similar to the step $n=4$ two preordering sequences ${\bf 5_1}\star {\bf 4_1}, {\bf 5_1}\star {\bf 3_1}\star {\bf 3_1}$ and ${\bf 6_2}\star{\bf 3_1}$ are omitted since these sequences are with knots which are repeats of the knots in the above preordering sequences. $\diamond$

Then to find the chains of transition for this step let us order these four preordering sequences with the following ordering where the knot ${\bf 3_1}\star {\bf 3_1}\star {\bf 3_1}\star {\bf 3_1}\star {\bf 3_1}$ assigned at the position of $2^5$ is put to follow the four sequences:
\begin{equation}
\begin{array}{rl}
& {\bf 6_1}\star {\bf 3_1}; \\
& {\bf 5_2}\star{\bf 4_1},\\
& {\bf 5_2}\star {\bf 3_1}\star {\bf 3_1};\\
& {\bf 4_1}\star{\bf 5_1}, \\
& {\bf 4_1}\star ({\bf 3_1}\star{\bf 4_1}),\\
& {\bf 4_1}\star {\bf 5_2}, \\
& {\bf 4_1}\star ({\bf 3_1}\star {\bf 3_1}\star {\bf 3_1});\\
& {\bf 3_1}\star ({\bf 3_1}\times {\bf 3_1}), \\
& {\bf 3_1}\star({\bf 3_1}\star {\bf 5_1}),\\
& {\bf 3_1}\star {\bf 6_1}, \\
& {\bf 3_1}\star ({\bf 3_1}\star{\bf 5_2}),\\
& {\bf 3_1}\star{\bf 6_2}, \\
& {\bf 3_1}\star ({\bf 4_1}\star {\bf 4_1}), \\
& {\bf 3_1}\star ({\bf 3_1}\star {\bf 3_1}\star {\bf 4_1}),\\
& {\bf 3_1}\star ({\bf 3_1}\star {\bf 3_1}\star {\bf 3_1}\star {\bf 3_1});\\
& {\bf 3_1}\star {\bf 3_1}\star {\bf 3_1}\star {\bf 3_1}\star {\bf 3_1}
\end{array}
\label{CC5}
\end{equation}

For this step we have three composite knots ${\bf 3_1}\star ({\bf 4_1}\star{\bf 5_1})$, ${\bf 5_1}\star {\bf 5_1}$ and ${\bf 4_1}\star ({\bf 4_1}\star{\bf 4_1})$ (related with $2\cdot3\cdot5$,$5\cdot5$ and $3\cdot3\cdot3$ respectively) of jumping over and there are two new knots ${\bf 4_1}\star {\bf 5_1}$ and ${\bf 3_1}\star ({\bf 3_1}\times{\bf 3_1})$ coming from the previous step. Thus there is a room for the introduction of new knot obtained only by the $\times$ operation. Then this new knot must be the composite knot ${\bf 3_1}\times {\bf 4_1}$ since besides the composite knot ${\bf 3_1}\times {\bf 3_1}$ it is the smallest of composite knots of this kind.

From (\ref{CC5})
there is a chain of transition given by $18\to 21\to 22\to 26\to 28 \to 27$ and the composite knot ${\bf 4_1}\star ({\bf 4_1}\star{\bf 4_1})$ related with $27=3\cdot3\cdot3$ is pushed out into the next step by the composite knot ${\bf 5_2}\star{\bf 4_1}$ at the starting position $18$. Then this repeated knot must be replaced by a new composite knot obtained by the $\times$ operation only and this new composite knot must be the knot ${\bf 3_1}\times {\bf 4_1}$. 

Then the composite knots at the intermediate states are assigned with the numbers of these states respectively.

In addition to the above chain there are two more chains: $24\to 30$ and $20\to 25$.
The
chain $24\to 30$ starts from ${\bf 3_1}\star ({\bf 3_1}\times{\bf 3_1})$ at $24$ and the composite knot ${\bf 3_1}\star ({\bf 4_1}\star{\bf 5_1})$ at $30$ is pushed out by the composite knot ${\bf 3_1}\star ({\bf 3_1}\star{\bf 3_1}\star{\bf 4_1})$. 
Then the chain $20\to 25$ starts from ${\bf 4_1}\star {\bf 5_1}$ at $20$ and  the composite knot 
${\bf 5_1}\star {\bf 5_1}$ at $25$ is pushed out by the composite knot ${\bf 3_1}\star ({\bf 3_1}\star{\bf 5_1})$.
 
Then the knots ${\bf 3_1}\star ({\bf 3_1}\star{\bf 3_1}\star{\bf 4_1})$ and ${\bf 3_1}\star ({\bf 3_1}\star{\bf 5_1})$ at the intemediate states of these two chains are assigned with the numbers $30=2\cdot3\cdot5$ and $25=5\cdot5$ respectively.

Now the remaining repeated composite knots at the positions $17,19, 23,29,31$ must be replaced by new prime knots  and thus $17,19, 23,29,31$ are determined to be prime numbers and they are determined by the prime numbers in the previous induction steps. Then we may follow the usual table of knots to determine that the new prime knots for the prime numbers $17,19, 23,29,31$ are ${\bf 6_3}$, ${\bf 7_1}$, ${\bf 7_2}$, ${\bf 7_3}$ and ${\bf 7_4}$ respectively (As the above induction steps we may just put in five prime knots and then later determine what these five prime knots will be. As shown in the above proof by using knots of the form ${\bf 3_1}\times{\bf 7_{(\cdot)}}$  we can determine the number of prime knots with minimal number of crossings $=7$. From this we can then determine these five prime knots). 

In summary we have the following form of the step $n=5$:
\begin{equation}
\begin{array}{rl}
& {\bf 6_3} \\
& {\bf 3_1}\times{\bf 4_1}\\
& {\bf 7_1}\\
& {\bf 4_1}\star{\bf 5_1} \\
& {\bf 4_1}\star ({\bf 3_1}\star{\bf 4_1})\\
& {\bf 4_1}\star {\bf 5_2} \\
& {\bf 7_2}\\
& {\bf 3_1}\star ({\bf 3_1}\times {\bf 3_1}) \\
& {\bf 3_1}\star({\bf 3_1}\star {\bf 5_1})\\
& {\bf 3_1}\star {\bf 6_1} \\
& {\bf 3_1}\star ({\bf 3_1}\star{\bf 5_2})\\
& {\bf 3_1}\star{\bf 6_2} \\
& {\bf 7_3} \\
& {\bf 3_1}\star ({\bf 3_1}\star {\bf 3_1}\star {\bf 4_1})\\
& {\bf 7_4}\\
& {\bf 3_1}\star {\bf 3_1}\star {\bf 3_1}\star {\bf 3_1}\star {\bf 3_1}
\end{array}
\label{CC6}
\end{equation}

This completes the induction step at $n=5$. We have that the structure of numbers of this step (including distribution of prime numbers in this step) is determined by the structure of numbers of the previous induction steps.
$\diamond$

\section{A Classification Table of Knots II}\label{sec2}

Following the above classification table up to $2^5$ let us in this section give the table up to $2^7$.
Again we shall see from the table that the preordering property is clear. At the $7$th step  there is a special composite knot ${\bf 4_1}\star{\bf 5_1}\star{\bf 5_1}$ which is not of jumping over and is not in the preordering sequences (On the other hand the knot ${\bf 5_1}\star{\bf 5_1}\star{\bf 5_1}$ is of jumping over).

We remark again that it is interesting that (by the ordering of composite knots with the $\times$ operation only) at the $6$th step we require exactly two prime knots with minimal number of crossings $=5$ to form the two composite knots obtained by the $\times$ operation only. From this we can determine the number of prime knots with minimal number of crossings $=5$ without using the actual contruction of these prime knots. We then denote these two prime knots by ${\bf 5_1}$ and ${\bf 5_2}$ respectively and the two composite knots obtained by the $\times$ operation only by ${\bf 3_1}\times{\bf 5_1}$ and ${\bf 3_1}\times{\bf 5_2}$ respectively. Similarly at the $7$th step we can determine that there are exactly three prime knots with minimal number of crossings $=6$ and we denote these three prime knots by ${\bf 6_1}$ and ${\bf 6_2}$ and ${\bf 6_3}$ respectively. These three prime knots give the composite knots ${\bf 3_1}\times{\bf 6_1}$, ${\bf 3_1}\times{\bf 6_2}$ and ${\bf 3_1}\times{\bf 6_3}$ respectively. We can then expect that at the next $8$th step we may determine that the number of prime knots with minimal number of crossings $=7$ is $7$ and then at the next $9$th step the number of prime knots with minimal number of crossings $=8$ is $21$, and so on; as we know from the well known table of prime knots \cite{Rol}.
Here the point is that we can determine the number of prime knots with the same minimal number of crossings without using the actual construction of these prime knots (and by using only the classification table of knots).
\begin{displaymath}
\begin{array}{|c|c|c|} \hline
\mbox{Type of Knot}& \mbox{Assigned number} \,\, |m|
 &\mbox{Repeated Knots being replaced}
\\ \hline

{\bf 3_1\star 6_3} & 33 & {\bf } \\ \hline

{\bf 3_1\star(3_1\times 4_1)} & 34 & {\bf } \\ \hline

{\bf 3_1\star 7_1} & 35 & {\bf } \\ \hline

{\bf 3_1\times 5_1} & 36 & {\bf 3_1\star(4_1\star 5_1)} \\ \hline

{\bf 7_5} & 37 & {\bf 3_1\star(4_1\star 3_1\star 4_1)} \\ \hline

{\bf 3_1\times 5_2} & 38 & {\bf 3_1\star(4_1\star 5_2)} \\ \hline

{\bf 3_1\star 7_2} & 39 & {\bf } \\ \hline

{\bf 3_1\star(3_1\star 3_1\times 3_1)} &40  & {\bf } \\ \hline

{\bf 7_6} & 41 & {\bf 3_1\star(3_1\star3_1\star 5_1)} \\ \hline

{\bf 5_1\star 5_1} & 42 & {\bf } \\ \hline

{\bf 7_7} & 43 & {\bf 5_1\star(3_1\star4_1)} \\ \hline

{\bf 5_1\star 5_2} & 44 & {\bf } \\ \hline

{\bf 4_1\times 4_1} & 45 & {\bf 5_1\star(3_1\star3_1\star3_1), 5_2\star(3_1\star4_1)} \\ \hline

{\bf 5_2\star5_2} & 46 & {\bf } \\ \hline

{\bf 8_1} & 47 & {\bf 5_2\star(3_1\star3_1\star3_1)},{\bf 4_1\star(3_1\times3_1)} \\ \hline

{\bf 4_1\star(3_1\star5_1)} & 48 & {\bf } \\ \hline

{\bf 4_1\star6_1} & 49 & {\bf } \\ \hline

{\bf 4_1\star(3_1\star5_2)} & 50 & {\bf } \\ \hline

{\bf 4_1\star6_2} & 51 & {\bf } \\ \hline

{\bf 4_1\star(4_1\star4_1)} & 52 & {\bf } \\ \hline

{\bf 8_2} & 53 & {\bf 4_1\star(4_1\star3_1\star3_1)} \\ \hline

{\bf 3_1\star(3_1\star3_1\star5_1)} & 54 & {\bf } \\ \hline

{\bf 3_1\star(3_1\star6_1)} & 55 & {\bf } \\ \hline

{\bf 3_1\star(3_1\star3_1\star5_2)} & 56 & {\bf } \\ \hline

{\bf 3_1\star(3_1\star6_2)} & 57 & {\bf } \\ \hline

{\bf 3_1\star7_3} & 58 & {\bf } \\ \hline

{\bf 8_3} & 59 & {\bf 3_1\star(3_1\star3_1\star3_1\star4_1)} \\ \hline

{\bf 3_1\star7_4} &60  & {\bf } \\ \hline

{\bf 8_4} &61  & {\bf 3_1\star(3_1\star3_1\star3_1\star 3_1\star3_1) } \\ \hline

{\bf 4_1\star(4_1\star3_1\star3_1)} & 62 & {\bf } \\ \hline

{\bf 4_1\star(3_1\star3_1\star3_1\star3_1)} & 63 & {\bf } \\ \hline

{\bf 3_1\star3_1\star3_1\star3_1\star3_1\star3_1} &64  & {\bf } \\ \hline
\end{array}
\end{displaymath}

\begin{displaymath}
\begin{array}{|c|c|c|} \hline
\mbox{Type of Knot}& \mbox{Assigned number} \,\, |m|
 &\mbox{Repeated Knots being replaced}
\\ \hline
{\bf 3_1\star(3_1\star6_3)} & 65 & {\bf } \\ \hline

{\bf 3_1\times(3_1\times3_1)} & 66 & {\bf 3_1\star(3_1\star3_1\times4_1)} \\ \hline

{\bf 8_5} & 67 & {\bf 3_1\star(3_1\star7_1)} \\ \hline

{\bf 4_1\times5_1} &68  & {\bf 3_1\star(3_1\star4_1\star5_1)} \\ \hline

{\bf 4_1\times(3_1\star4_1)} &69  & {\bf 3_1\star(3_1\star4_1\star3_1\star4_1)} \\ \hline

{\bf 4_1\times5_2} & 70 & {\bf 3_1\star(3_1\star4_1\star5_2)} \\ \hline

{\bf 8_6} & 71 & {\bf 3_1\star(3_1\star7_2)} \\ \hline

{\bf 4_1\star6_3} & 72 & {\bf } \\ \hline

{\bf 8_7} &73  & {\bf 4_1\star(3_1\times4_1)} \\ \hline

{\bf 5_1\star(3_1\times3_1)} &74  & {\bf } \\ \hline

{\bf 5_1\star(3_1\star5_1)} & 75 & {\bf } \\ \hline

{\bf 5_1\star6_1} & 76 & {\bf } \\ \hline

{\bf 5_1\star(3_1\star5_2)} & 77 & {\bf } \\ \hline

{\bf 5_1\star6_2} &78  & {\bf } \\ \hline

{\bf 8_8} & 79 & {\bf 5_1\star(4_1\star4_1), 5_2\star(3_1\star5_1)} \\ \hline

{\bf 5_2\star6_1} &80  & {\bf } \\ \hline

{\bf 5_2\star(3_1\star5_2)} & 81 & {\bf } \\ \hline

{\bf 5_2\star6_2} &82  & {\bf } \\ \hline

{\bf 8_9} &83  & {\bf 5_2\star(4_1\star4_1)} \\ \hline

{\bf 4_1\star7_1} &84  & {\bf } \\ \hline

{\bf 4_1\star(4_1\star5_1)} & 85 & {\bf } \\ \hline

{\bf 4_1\star(4_1\star4_1\star3_1)} &86  & {\bf } \\ \hline

{\bf 4_1\star(4_1\star5_2)} & 87 & {\bf } \\ \hline

{\bf 4_1\star7_2} &88  & {\bf } \\ \hline

{\bf 8_{10}} & 89 & {\bf 4_1\star(3_1\star3_1\times3_1)} \\ \hline

{\bf 4_1\star(3_1\star3_1\star5_1)} &90  & {\bf } \\ \hline

{\bf 4_1\star(3_1\star6_1)} & 91 & {\bf } \\ \hline

{\bf 4_1\star(3_1\star3_1\star5_2)} & 92 & {\bf } \\ \hline

{\bf 4_1\star(3_1\star6_2)} & 93 & {\bf } \\ \hline

{\bf 4_1\star7_3} & 94 & {\bf } \\ \hline

{\bf 4_1\star(4_1\star3_1\star3_1\star3_1)} &95  & {\bf } \\ \hline

{\bf 4_1\star7_4} & 96 & {\bf } \\ \hline
\end{array}
\end{displaymath}

\begin{displaymath}
\begin{array}{|c|c|c|} \hline
\mbox{Type of Knot}& \mbox{Assigned number}\,\, |m|
 &\mbox{Repeated Knots being replaced}
\\ \hline

{\bf 8_{11}} & 97 & {\bf 4_1\star(3_1\star3_1\star3_1\star3_1\star3_1)} \\ \hline

{\bf 4_1\star(5_1\star5_1)} &98  & {\bf 3_1\star(3_1\star6_3)} \\ \hline

{\bf 3_1\star(3_1\star3_1\times4_1)} &99  & {\bf } \\ \hline

{\bf 3_1\star(3_1\star7_1)} &100  & {\bf } \\ \hline

{\bf 8_{12}} & 101 & {\bf3_1\star(3_1\times5_1) } \\ \hline

{\bf 3_1\star7_5} & 102 & {\bf } \\ \hline

{\bf 8_{13}} &103  & {\bf 3_1\star(3_1\times5_2)} \\ \hline

{\bf 3_1\star(3_1\star7_2)} & 104 & {\bf } \\ \hline

{\bf 3_1\star(3_1\star3_1\star3_1\times3_1)} & 105 & {\bf } \\ \hline

{\bf 3_1\star7_6} & 106 & {\bf } \\ \hline

{\bf 8_{14}} &107  & {\bf 3_1\star(5_1\star5_1)} \\ \hline

{\bf 3_1\star7_7} & 108 & {\bf } \\ \hline

{\bf 8_{15}} &109  & {\bf 3_1\star(5_1\star5_2)} \\ \hline

{\bf 3_1\times6_1} &110  & {\bf 3_1\star(5_1\star5_2)} \\ \hline

{\bf 3_1\times(3_1\star5_2)} & 111 & {\bf 3_1\star(4_1\times4_1)} \\ \hline

{\bf 3_1\times6_2} &112 & {\bf 3_1\star(5_2\star5_2)} \\ \hline

{\bf 8_{16}} & 113 & {\bf 3_1\star(5_2\star5_2)} \\ \hline

{\bf 3_1\star8_1} & 114 & {\bf } \\ \hline

{\bf 3_1\times6_3} & 115 & {\bf 3_1\star(4_1\star3_1\star5_1),3_1\star(4_1\star4_1\star4_1)} \\ \hline

{\bf 3_1\star8_2} &116  & {\bf } \\ \hline

{\bf 3_1\star(3_1\star3_1\star3_1\star5_1)} &117  & {\bf } \\ \hline

{\bf 3_1\star(3_1\star3_1\star6_1)} & 118 & {\bf } \\ \hline

{\bf 3_1\star(3_1\star3_1\star3_1\star5_2)} & 119 & {\bf } \\ \hline

{\bf 3_1\star(3_1\star3_1\star6_2)} &120  & {\bf } \\ \hline

{\bf 3_1\star(3_1\star3_1\star7_3)} &121  & {\bf } \\ \hline

{\bf 3_1\star8_3} & 122 & {\bf } \\ \hline

{\bf 3_1\star(3_1\star7_4)} &123  & {\bf } \\ \hline

{\bf 3_1\star8_4} &124  & {\bf } \\ \hline

{\bf 3_1\times(3_1\times4_1)} &125  & {\bf 3_1\star(4_1\star4_1\star3_1\star3_1)} \\ \hline

{\bf 3_1\star(4_1\star3_1\star3_1\star3_1\star3_1)} & 126 & {\bf } \\ \hline

{\bf 8_{17}} &127  & {\bf 3_1\star(3_1\star3_1\star3_1\star3_1\star3_1\star3_1)} \\ \hline

{\bf 3_1\star(3_1\star3_1\star3_1\star3_1\star3_1\star3_1)} & 128 & {\bf } \\ \hline

\end{array}
\end{displaymath}

\section{Applications in Number Theory}\label{sec3}

Let us apply the results in the above sections to solve some problems in number theory. Let us prove the following Goldbach Conjecture \cite{Har}-\cite{Che3}:
\begin{theorem}(Goldbach Conjecture)

Each even number $\geq 6$ can be written as the sum of two prime numbers.
\end{theorem}

{\bf Remark}. The statement of the Goldbach Conjecture gives a correlation of numbers and prime numbers. Let us call a property of numbers and prime numbers as a correlated property if it gives a correlation of numbers and prime numbers. Then the statement of the Goldbach Conjecture whenever true is a correlated property of numbers and prime numbers.

{\bf Proof of the Goldbach Conjecture}. 
Let us first consider a general correlated property of numbers and prime numbers.
In the above sections and the above theorem \ref{maintheorem} we have proved that 
the structure of numbers (obtained by assigning knots to numbers) of the coming induction steps (including the distribution of prime numbers in the coming induction steps) 
is by induction determined by the structure of numbers (obtained by assigning knots to numbers) of the previous steps. This structure of numbers is the original number system together with the one-to-one assignment of numbers to knots.
 
 Further as shown in the proof of theorem \ref{maintheorem} the structure of numbers of the previous steps are inherited in the structure of numbers of the coming induction steps as the preordering sequences in the determination of the structure of numbers of the coming induction steps. Thus all the properties of the structure of numbers of the previous steps are inherited in the structure of numbers of the coming induction steps in the determination of the structure of numbers of the coming induction steps.
 
It follows that the structure of numbers of the coming induction steps inherits all the properties of the structure of numbers of the previous steps.

 Now since the correlated properties of numbers and prime numbers in the previous steps are properties of the structure of numbers of the previous steps
we have that the structure of numbers of the coming induction steps inherits all the correlated properties of numbers and prime numbers of the previous steps.

Thus the correlated properties of numbers and prime numbers in the previous steps must be extended as the same correlated properties for the numbers and prime numbers in the coming induction steps. 

Thus the forms of these extended correlated properties for the numbers and prime numbers in the coming induction steps must be analogous to the forms of the corresponding correlated properties of numbers and prime numbers in the previous  steps and are extended from the forms of the corresponding correlated properties of numbers and prime numbers in the previous steps.

Further since each induction step is mixed with new prime numbers and composite numbers composed with prime numbers in the previous steps we have that the composite numbers and prime numbers in the previous 
 steps are correlated with the composite numbers and prime numbers in the coming induction steps 
 (where the correlation is from the correlated properties of composite numbers and prime numbers in the previous steps). 
 Thus we have that the extended correlated properties for the composite numbers and prime numbers in the coming induction steps must overlap with the corresponding correlated properties of composite numbers and prime numbers in the previous
steps. This form of overlapping then  gives a guide to find out the extended correlated properties for the composite numbers and prime numbers in the coming induction steps.

Now let us consider the case of the Goldbach Conjecture.

At the step $n=3$ we have the following correlated property of numbers and prime numbers that each even number $\geq 6$ in the steps $n=1,2,3$ can be written as the sum of two prime numbers (in the steps $n=1,2,3$) which are from a pair of twin prime numbers (in the steps $n=1,2,3$) respectively:
\begin{equation}
6=3+3, 8=3+5
\label{gold1}
\end{equation}
where $3$ and $5$ is a pair of twin prime numbers and $5$ and $7$ is another pair of twin prime numbers and that there is a pair of twin prime numbers $5$ and $7$ in the step $n=3$ (and not in the steps $n=1,2$).

Further there are even numbers in the next induction step $n=4$ which can be written as the sum of two prime numbers (in the steps $n=1,2,3$) which are from a pair of twin prime numbers (in the steps $n=1,2,3$) respectively:
\begin{equation}
10=5+5, 12=5+7, 14=7+7
\label{gold2}
\end{equation}
where $10, 12,14$ are even numbers in the next induction step $n=4$ and $5$ and $7$ is a pair of twin prime numbers in the step $n=3$ (and not in the step $n=1,2$). Thus (\ref{gold2}) is the overlapping of the step $n=3$ and the induction step $n=4$. This correlated property of prime numbers in the steps $n=1,2,3$ is as the Goldbach Conjecture at the step $n=3$. In summary we have the following form of the Goldbach Conjecture at the step $n=3$ (which is already a correlated property of numbers and prime numbers in the steps $n=3$):

{\bf Statement of the Goldbach Conjecture at the step ${\bf n=3}$}. Each even number in the step $n=3$ can be written as the sum of two prime numbers in the steps $n=1,2,3$
and these two prime numbers are each from a pair of twin prime numbers in the steps $n=1,2,3$ and that there is a pair of twin prime numbers in the step $n=3$ (and not in the step $n=1,2$) such that the sum of this pair of twin prime numbers is an even number in the step $n=4$.

Further it is clear that the sum of two prime numbers which are each from a pair of twin prime numbers  in the step $n=3$ (and not in the step $n=1,2$) is an even number in the step $n=4$. This is as the overlapping of the step $n=3$ and step $n=4$. $\diamond$

Now we consider the following form of the Goldbach Conjecture at the step $n=4$:

{\bf Statement of the Goldbach Conjecture at the step ${\bf n=4}$}. Each even number in the step $n=4$ can be written as the sum of two prime numbers in the steps $n=1,2,3,4$
and these two prime numbers are each from a pair of twin prime numbers in the steps $n=1,2,3,4$ and that there is a pair of twin prime numbers in the step $n=4$ (and not in the steps $n=1,2,3$) such that the sum of this pair of twin prime numbers is an even number in the step $n=5$.

Further it is clear that the sum of two prime numbers which are each from a pair of twin prime numbers in the step $n=4$ (and not in the steps $n=1,2,3$) is an even number in the step $n=5$. This is as the overlapping of the step $n=4$ and step $n=5$. $\diamond$

From the overlapping in the step $n=3$ we have that the extension of the statement of this overlapping to the step $n=4$ is that the even numbers in the step $n=4$ can be written as the sum of two prime numbers which are each from a pair of twin prime numbers in the step $n=4$. 
Then this statement includes the statement that there is a pair of twin prime numbers in the step $n=4$ (and not in the steps $n=1,2,3$) which is in the statement of the Goldbach Conjecture in the step $n=4$ since there must exist an even muber (e.g. the even number $2^4$) in the step $n=4$ which cannot be written as the sum of prime numbers in the steps $n=1,2,3$.

Thus together with step $n=3$ the extension of the statement of the overlapping to the step $n=4$ is that
each even number in the step $n=4$ can be written as the sum of two prime numbers and these two prime numbers are each from a pair of twin prime numbers in the steps $n=1,2,3,4$ and that there is a pair of twin prime numbers in the step $n=4$ (and not in the step $n=1,2,3$) such that the sum of this pair of twin prime numbers is an even number in the step $n=5$.
This is just the statement of the Goldbach Conjecture in the step $n=4$ since the statement of the overlapping in the statement of the Goldbach Conjecture in the step $n=4$ is always true. This shows that the statement of the Goldbach Conjecture at the step $n=4$ is the unique extension to the step $n=4$ of the overlapping of the step $n=3$ and the step $n=4$.

It follows that the statement of the Goldbach Conjecture at the step $n=4$ is  just the correlated property of the numbers and prime numbers in the step $n=4$ extended from the 
Goldbach Conjecture at the step $n=3$ which is a correlated property of the numbers and prime numbers in the steps $n=1,2,3$. This proves that the Goldbach Conjecture at the step $n=4$ holds.

 Then by the same induction procedure we have that the statement of the Goldbach Conjecture 
 at the induction step $n=5$ (which is analogous to step $n=3$ and step $n=4$) must hold. Continuing in this way we have that the Goldbach Conjecture must hold at each $n$th step. Thus we have proved that each even number $\geq 6$ can be written as the sum of two prime numbers which are from a pair of twin prime numbers respectively and that there is a pair of twin prime numbers in the $n$th step (and not in the previous steps). This statement implies the usual Goldbach Conjecture and thus the Goldbach Conjecture holds, as was to be proved. $\diamond$

From this proof of the Goldbach Conjecture we can also prove the following Twin Prime Conjecture:
\begin{theorem}(Twin Prime Conjecture)
There exist infinitely many pairs of twin prime numbers.
\end{theorem}
{\bf Proof}. In the above proof we have proved that there is a pair of twin prime numbers in the $n$th step (and not in the previous steps). Thus there exist infinitely many pair of twin prime numbers, as was to be proved. $\diamond$

\section{Conclusion}\label{sec16}

By using the connected sum operations $\star$ and $\times$ on knots a classification table of knots can be formed such that knots are one-to-one assigned with an integer and prime knots are bijectively assigned with prime numbers such that the prime number $2$ corresponds to the trefoil knot. This assignment is onto the set of positive integers except $2$ where the trefoil knot is assigned with 1 and is related to $2$. Further this assignment for the $n$th induction step of the number $2^n$ is determined by this assignment for the previous $n-1$ steps. From this induction of assigning knots with numbers we can prove the Goldbach Conjecture and the Twin Prime Conjecture.  
This induction of assigning knots with numbers may also be used to investigate other problems in number theory.

\end{document}